\documentclass[11pt]{amsart}    % [a4paper, 11pt]
\usepackage{latexsym, amsmath, amsthm, amssymb, setspace, verbatim}
\usepackage[all]{xy}
\usepackage{longtable}
\usepackage[applemac]{inputenc}
\usepackage{hyperref}
\usepackage{color}
\usepackage[usenames,dvipsnames]{xcolor}

\title[ The continuous Rokhlin property and the UCT ]{ A short note on the continuous Rokhlin property and the universal coefficient theorem}
\author{ Gábor Szabó }
\address{Westfälische Wilhelms-Universität, Fachbereich Mathematik, \phantom{--------------}\linebreak \text{}\hspace{3.5mm} Einsteinstrasse 62, 48149 Münster, Germany}
\email{gabor.szabo@uni-muenster.de}
\thanks{\emph{Supported by:} SFB 878 \emph{Groups, Geometry and Actions} }

\begin{document}

% math
\renewcommand\matrix[1]{\left(\begin{array}{*{10}{c}} #1 \end{array}\right)}  % Matrix
\newcommand\set[1]{\left\{#1\right\}}  % Menge
\newcommand\mset[1]{\left\{\!\!\left\{#1\right\}\!\!\right\}}

%% Besondere Variablen
%Zahlmengen-Stil
\newcommand{\IA}[0]{\mathbb{A}} \newcommand{\IB}[0]{\mathbb{B}}
\newcommand{\IC}[0]{\mathbb{C}} \newcommand{\ID}[0]{\mathbb{D}}
\newcommand{\IE}[0]{\mathbb{E}} \newcommand{\IF}[0]{\mathbb{F}}
\newcommand{\IG}[0]{\mathbb{G}} \newcommand{\IH}[0]{\mathbb{H}}
\newcommand{\II}[0]{\mathbb{I}} \renewcommand{\IJ}[0]{\mathbb{J}}
\newcommand{\IK}[0]{\mathbb{K}} \newcommand{\IL}[0]{\mathbb{L}}
\newcommand{\IM}[0]{\mathbb{M}} \newcommand{\IN}[0]{\mathbb{N}}
\newcommand{\IO}[0]{\mathbb{O}} \newcommand{\IP}[0]{\mathbb{P}}
\newcommand{\IQ}[0]{\mathbb{Q}} \newcommand{\IR}[0]{\mathbb{R}}
\newcommand{\IS}[0]{\mathbb{S}} \newcommand{\IT}[0]{\mathbb{T}}
\newcommand{\IU}[0]{\mathbb{U}} \newcommand{\IV}[0]{\mathbb{V}}
\newcommand{\IW}[0]{\mathbb{W}} \newcommand{\IX}[0]{\mathbb{X}}
\newcommand{\IY}[0]{\mathbb{Y}} \newcommand{\IZ}[0]{\mathbb{Z}}

%Geschwungener Stil
\newcommand{\CA}[0]{\mathcal{A}} \newcommand{\CB}[0]{\mathcal{B}}
\newcommand{\CC}[0]{\mathcal{C}} \newcommand{\CD}[0]{\mathcal{D}}
\newcommand{\CE}[0]{\mathcal{E}} \newcommand{\CF}[0]{\mathcal{F}}
\newcommand{\CG}[0]{\mathcal{G}} \newcommand{\CH}[0]{\mathcal{H}}
\newcommand{\CI}[0]{\mathcal{I}} \newcommand{\CJ}[0]{\mathcal{J}}
\newcommand{\CK}[0]{\mathcal{K}} \newcommand{\CL}[0]{\mathcal{L}}
\newcommand{\CM}[0]{\mathcal{M}} \newcommand{\CN}[0]{\mathcal{N}}
\newcommand{\CO}[0]{\mathcal{O}} \newcommand{\CP}[0]{\mathcal{P}}
\newcommand{\CQ}[0]{\mathcal{Q}} \newcommand{\CR}[0]{\mathcal{R}}
\newcommand{\CS}[0]{\mathcal{S}} \newcommand{\CT}[0]{\mathcal{T}}
\newcommand{\CU}[0]{\mathcal{U}} \newcommand{\CV}[0]{\mathcal{V}}
\newcommand{\CW}[0]{\mathcal{W}} \newcommand{\CX}[0]{\mathcal{X}}
\newcommand{\CY}[0]{\mathcal{Y}} \newcommand{\CZ}[0]{\mathcal{Z}}

%Script Stil
\newcommand{\FA}[0]{\mathfrak{A}} \newcommand{\FB}[0]{\mathfrak{B}}
\newcommand{\FC}[0]{\mathfrak{C}} \newcommand{\FD}[0]{\mathfrak{D}}
\newcommand{\FE}[0]{\mathfrak{E}} \newcommand{\FF}[0]{\mathfrak{F}}
\newcommand{\FG}[0]{\mathfrak{G}} \newcommand{\FH}[0]{\mathfrak{H}}
\newcommand{\FI}[0]{\mathfrak{I}} \newcommand{\FJ}[0]{\mathfrak{J}}
\newcommand{\FK}[0]{\mathfrak{K}} \newcommand{\FL}[0]{\mathfrak{L}}
\newcommand{\FM}[0]{\mathfrak{M}} \newcommand{\FN}[0]{\mathfrak{N}}
\newcommand{\FO}[0]{\mathfrak{O}} \newcommand{\FP}[0]{\mathfrak{P}}
\newcommand{\FQ}[0]{\mathfrak{Q}} \newcommand{\FR}[0]{\mathfrak{R}}
\newcommand{\FS}[0]{\mathfrak{S}} \newcommand{\FT}[0]{\mathfrak{T}}
\newcommand{\FU}[0]{\mathfrak{U}} \newcommand{\FV}[0]{\mathfrak{V}}
\newcommand{\FW}[0]{\mathfrak{W}} \newcommand{\FX}[0]{\mathfrak{X}}
\newcommand{\FY}[0]{\mathfrak{Y}} \newcommand{\FZ}[0]{\mathfrak{Z}}

\newcommand{\fc}[0]{\mathfrak{c}}

%Pfeilbefehle abkürzen
\newcommand{\Ra}[0]{\Rightarrow}
\newcommand{\La}[0]{\Leftarrow}
\newcommand{\LRa}[0]{\Leftrightarrow}

%Modifikation der Variablen
\renewcommand{\phi}[0]{\varphi}
\newcommand{\eps}[0]{\varepsilon}

%zusätzliche Features
\newcommand{\quer}[0]{\overline}
\newcommand{\uber}[0]{\choose}
\newcommand{\ord}[0]{\operatorname{ord}}		% Ordnung
\newcommand{\GL}[0]{\operatorname{GL}}
\newcommand{\supp}[0]{\operatorname{supp}}	% Träger
\newcommand{\id}[0]{\operatorname{id}}		% Identität
\newcommand{\Sp}[0]{\operatorname{Sp}}		% Spektrum eines Elements
\newcommand{\eins}[0]{\mathbf{1}}			% Eine Eins in allgemeinerem Kontext, z.B. in einem Ring
\newcommand{\diag}[0]{\operatorname{diag}}
\newcommand{\auf}[1]{\quad\stackrel{#1}{\longrightarrow}\quad}
\newcommand{\prim}[0]{\operatorname{Prim}}
\newcommand{\ad}[0]{\operatorname{Ad}}
\newcommand{\ext}[0]{\operatorname{Ext}}
\newcommand{\ev}[0]{\operatorname{ev}}
\newcommand{\fin}[0]{{\subset\!\!\!\subset}}
\newcommand{\diam}[0]{\operatorname{diam}}
\newcommand{\Hom}[0]{\operatorname{Hom}}
\newcommand{\Aut}[0]{\operatorname{Aut}}
\newcommand{\del}[0]{\partial}
\newcommand{\dimnuc}[0]{\dim_{\mathrm{nuc}}}
\newcommand{\dr}[0]{\operatorname{dr}}
\newcommand{\dimrok}[0]{\dim_{\mathrm{Rok}}}
\newcommand{\dimrokcyc}[0]{\dim_{\mathrm{Rok}}^{\mathrm{cyc}}}
\newcommand{\dimrokcycc}[0]{\dim_{\mathrm{Rok}}^{\mathrm{cyc,c}}}
\newcommand{\dimnuceins}[0]{\dimnuc^{\!+1}}
\newcommand{\dreins}[0]{\dr^{\!+1}}
\newcommand{\dimrokeins}[0]{\dimrok^{\!+1}}
\newcommand{\reldimrok}[2]{\dimrok(#1~|~#2)}
\newcommand{\mdim}[0]{\operatorname{mdim}}
\newcommand*\onto{\ensuremath{\joinrel\relbar\joinrel\twoheadrightarrow}} % surjectiver Pfeil
\newcommand*\into{\ensuremath{\lhook\joinrel\relbar\joinrel\rightarrow}}  % injektiver Pfeil
\newcommand{\im}[0]{\operatorname{Im}}
\newcommand{\dst}[0]{\displaystyle}
\newcommand{\cstar}[0]{$\mathrm{C}^*$}
\newcommand{\dist}[0]{\operatorname{dist}}
\newcommand{\ue}[0]{{\approx_{\mathrm{u}}}}
\newcommand{\mue}[0]{\approx_{\mathrm{mu}}}
\newcommand{\End}[0]{\operatorname{End}}
\newcommand{\Ell}[0]{\operatorname{Ell}}
\newcommand{\gr}[0]{\operatorname{gr}}
\newcommand{\Ost}[0]{\mathcal{O}_\infty^{\mathrm{st}}}
\newcommand{\Bst}[0]{\mathcal{B}^{\mathrm{st}}}
\newcommand{\inv}[0]{\operatorname{Inv}}
\newcommand{\ann}[0]{\operatorname{Ann}}

% theorems
\newtheorem{satz}{Satz}[section]		% <--- optional, zählt so mit den Abschnitten
\newtheorem{cor}[satz]{Corollary}
\newtheorem{lemma}[satz]{Lemma}
\newtheorem{prop}[satz]{Proposition}
\newtheorem{theorem}[satz]{Theorem}
\newtheorem*{theoreme}{Theorem}

\theoremstyle{definition}
\newtheorem{defi}[satz]{Definition}
\newtheorem*{defie}{Definition}
\newtheorem{defprop}[satz]{Definition \& Proposition}
\newtheorem{nota}[satz]{Notation}
\newtheorem*{notae}{Notation}
\newtheorem{rem}[satz]{Remark}
\newtheorem*{reme}{Remark}
\newtheorem{example}[satz]{Example}
\newtheorem{defnot}[satz]{Definition \& Notation}
\newtheorem{question}[satz]{Question}
\newtheorem*{questione}{Question}

\newenvironment{bew}{\begin{proof}[Proof]}{\end{proof}}

\begin{abstract} 
Let $G$ be a metrizable compact group, $A$ a separable \cstar-algebra and $\alpha: G\to\Aut(A)$ a strongly continuous action. Provided that $\alpha$ satisfies the continuous Rokhlin property, we show that the property of satisfying the UCT passes from $A$ to the crossed product \cstar-algebra $A\rtimes_\alpha G$ and the fixed point algebra $A^\alpha$. This extends a result by Gardella in the case that $G$ is the circle and $A$ is nuclear. For circle actions on separable, unital \cstar-algebras with the continuous Rokhlin property, we establish a connection between the $E$-theory equivalence class of the coefficient algebra $A$ and the fixed point algebra $A^\alpha$.
\end{abstract}

\maketitle

%\thispagestyle{empty}
%\newpage \tableofcontents
%\setcounter{page}{1}

\setcounter{section}{-1}

\section{Introduction}
\noindent
Within the field of \cstar-dynamical systems, the discovery and systematic study of various kinds of Rokhlin-type properties has recently become the driving force behind many recent and interesting results. For a loose selection of such instances, see \cite{ Izumi04, Izumi04II, Nawata13, Santiago14, Phillips11, Phillips12, OsakaPhillips12, Sato10, MatuiSato12, MatuiSato14, BarlakSzabo14, HirshbergWinterZacharias14, BarlakEndersMatuiSzaboWinter, Szabo14}.
Particularly notable for the purpose of this note is the investigation of compact group actions with the Rokhlin property (see \cite{HirshbergWinter07}) and more specially circle actions with the Rokhlin property, see \cite{Gardella14_1, Gardella14_2}. 

Despite the fact that many known results from the realm of finite group actions with the Rokhlin property carry over to the setting of compact groups, there remain some subtle difficulties concerning certain properties like the permanence of the UCT. For example, it is known that in many cases, the UCT passes from a \cstar-algebra to its crossed product associated to a Rokhlin action of a finite group, see \cite[Theorem 3]{Santiago14}. At present, it is unknown whether the UCT passes in all cases, or if such a permanence property holds for actions of compact groups. For example, when confronted with the problem of classifying Rokhlin actions of the circle on Kirchberg algebras by means of $K$-theory, this proves to be a rather annoying obstacle. To overcome this, Gardella introduced the continuous Rokhlin property for compact group actions on unital \cstar-algebras in \cite{Gardella14_2}. It turns out that this stronger version is compatible with the UCT for circle actions on nuclear \cstar-algebras. 
In this short note, we extend the definition of the continuous Rokhlin property to cover the non-unital case.
We then prove $E$-theoretic versions of Gardella's UCT preservation theorem for actions of all metrizable compact groups on separable \cstar-algebras by using a somewhat more conceptual approach, enabling less complicated proofs. At last, we show that for any strongly continuous \cstar-dynamical system $(A,\alpha,\IT)$ with the continuous Rokhlin property on a unital \cstar-algebra, we have that $A$ is $E$-equivalent to $A^\alpha\oplus SA^\alpha$. This particularly yields an $E$-theoretic version of Gardella's observation that $K_*(A)\cong K_*(A^\alpha)\oplus K_{*+1}(A^\alpha)$, whenever $\alpha: \IT\to\Aut(A)$ has the ordinary Rokhlin property.

I would like to express my gratitude to Sel\c{c}uk Barlak, who contributed some improvements to this note.

%%%%%%%%%%%%%%%%%%%%%%%%%%%%%%%%%%%%%%%%%%%%%%%%%

\section{The continuous Rokhlin property}

\defi Let $A$ be a separable \cstar-algebra. We denote the path algebra of $A$ by
\[
A_\fc = \CC_b\bigl( [1,\infty), A \bigl)/\CC_0\bigl( [1,\infty), A \bigl).
\]
Let $G$ be a metrizable, locally compact group and $\alpha: G\to\Aut(A)$ a strongly continuous action. Let $\CC_{b,\alpha}\bigl( [1,\infty), A \bigl)$ be the \cstar-subalgebra of $\CC_{b}\bigl( [1,\infty), A \bigl)$ consisting of those points on which the induced action of $\alpha$ is continuous. One sees easily that $\CC_{0}\bigl( [1,\infty), A \bigl)\subset \CC_{b,\alpha}\bigl( [1,\infty), A \bigl)$. We call
\[
A_{\fc,\alpha} = \CC_{b,\alpha}\bigl( [1,\infty), A \bigl)/\CC_{0}\bigl( [1,\infty), A \bigl)
\]
the continuous path algebra of $A$ with respect to $\alpha$. Clearly, $A$ always embeds as constant paths both into $A_\fc$ and $A_{\fc,\alpha}$. Similarly as in the case of sequence algebras (see \cite{Kirchberg04}), we define the central path algebra of $A$ by
\[
F_{\fc}(A) = A_\fc\cap A'/\ann(A,A_\fc)
\]
and the continuous central path algebra of $A$ with respect to $\alpha$
\[
F_{\fc,\alpha}(A) = A_{\fc,\alpha}\cap A'/\ann(A,A_{\fc,\alpha}).
\]

\rem If $G$ is compact, then the \cstar-algebras $F_\fc(A)$ and $F_{\fc,\alpha}(A)$ are unital. Let $h\in A^\alpha$ be a positive contraction that is strictly positive in $A$. Then the class of the path $b_t = h^{1/t}$ for $t\geq 1$ defines a unit for both $F_\fc(A)$ and $F_{\fc,\alpha}(A)$. 

\rem \label{central tensor}
In the above setting, we get two well-defined $*$-homomorphisms via
\[
A\otimes_{\max} F_\fc(A)\to A_\fc,\quad a\otimes\bigl( x+\ann(A,A_\fc) \bigl) ~\mapsto ax
\]
and
\[
A\otimes_{\max} F_{\fc,\alpha}(A)\to A_{\fc,\alpha},\quad a\otimes\bigl( x+\ann(A,A_{\fc,\alpha}) \bigl) ~\mapsto ax.
\]
Moreover, the image of $a\otimes\eins$ is $a$ under both these maps.

\rem \label{path action}
By definition, the action $\alpha$ on $A$ clearly extends to the action $\alpha_\fc$ on $A_{\fc,\alpha}$ given by 
\[
\alpha_{\fc,g}([(b_t)_{t\geq 1}]) = [\bigl( \alpha_g(b_t) \bigl)_{t\geq 1}]\quad\text{for all}~g\in G~\text{and}~[(b_t)_{t\geq 1}]\in A_{\fc,\alpha}.
\]
Since one easily checks that both $A_{\fc,\alpha}\cap A'$ and $\ann(A,A_{\fc,\alpha})$ are $\alpha_\fc$-invariant \cstar-subalgebras, we also get an induced action $\tilde{\alpha}_\fc$ on $F_{\fc,\alpha}(A)$ via
\[
\tilde{\alpha}_{\fc}\bigl( x+\ann(A,A_{\fc,\alpha}) \bigl) = \alpha_\fc(x)+\ann(A,A_{\fc,\alpha}).
\]
Having made these definitions, one checks that the above $*$-homomorphism from $A\otimes_{\max} F_{\fc,\alpha}(A)$ is $\alpha\otimes\tilde{\alpha}_{\fc}$-to-$\alpha_{\fc}$-equivariant.

\nota Given a compact group $G$, we denote by $\sigma$ the canonical $G$-shift action on $\CC(G)$ given by $\sigma_g(f)(h) = f(g^{-1}h)$ for all $g,h\in G$ and $f\in\CC(G)$.

\defi Let $A$ be a separable \cstar-algebra and $G$ a compact group. Let $\alpha: G\to\Aut(A)$ be a strongly continuous action. We say that $\alpha$ has the continuous Rokhlin property, if there exists a unital and equivariant $*$-homomorphism
\[
\phi: (\CC(G),\sigma) \to \bigl( F_{\fc,\alpha}(A), \tilde{\alpha}_\fc \bigl).
\]

\rem Let $(A,\alpha,G)$ be a strongly continuous \cstar-dynamical system as above. In the special case that $G=\IT$, the action $\alpha$ has the continuous Rokhlin property if and only if there is a unitary $u\in F_{\fc,\alpha}(A)$ with $\tilde{\alpha}_{\fc,\xi}(u)=\xi\cdot u$ for all $\xi\in\IT$. This follows from the observation that $u=\phi(\id_{\CC(\IT)})$ yields such a unitary, and conversely every such unitary generates an equivariant copy of the dynamical system $(\CC(\IT),\sigma)$.
In particular, this explains the analogy of the above definition to \cite[3.1]{Gardella14_2}.

%%%%%%%%%%%%%%%%%%%%%%%%%%

\section{Preservation of the UCT}

In this section, we show that the UCT is preserved under forming crossed products or fixed point algebras associated to compact group actions with the continuous Rokhlin property. We remark that within the context of this note, we exclusively work with the $E$-theory version of the UCT. 

\defi[compare to {\cite[Section 7]{GuentnerHigsonTrout}}] 
A separable \cstar-algebra $A$ satisfies the UCT in $E$-theory, if for every separable \cstar-algebra $B$, the natural map $E_*(A,B) \to \Hom\bigl( K_*(A), K_*(B) \bigl)$ gives rise to a short exact sequence
\[
\xymatrix@C-2mm{
0 \ar[r] & \operatorname{Ext}_\IZ^1\bigl( K_*(A), K_*(B) \bigl) \ar[r] & E_*(A,B) \ar[r] & \Hom\bigl( K_*(A), K_*(B) \bigl) \ar[r] & 0.
}
\]
As is the case for $KK$-theory, this is (a posteriori) the same as being $E$-equivalent to an abelian \cstar-algebra. For nuclear \cstar-algebras, $E$-theory coincides with $KK$-theory, so this makes no difference to the $KK$-theoretic version of the UCT.\vspace{2mm}

The following result is proved in \cite[23.10.8]{BlaKK} for $KK$-theory, but an analogous proof yields the same for $E$-theory.

{
\prop \label{E dominance}
Let $A$ and $B$ be separable \cstar-algebras. Assume that $A$ $E$-dominates $B$, i.e. there are $x\in E(B,A)$ and $y\in E(A,B)$ with $x\otimes y=1$ in $E(B,B)$. If $A$ satisfies the UCT, then so does $B$. 
}

\nota Let $(A,\alpha,G)$ and $(B,\beta,G)$ be two strongly continuous \cstar-dynamical systems. For an equivariant $*$-homomorphism $\psi: (A,\alpha)\to (B,\beta)$, we denote by $\hat{\psi}: A\rtimes_\alpha G\to B\rtimes_\beta G$ the induced $*$-homomorphism between the crossed products.

{
\theorem \label{cRok crossed product UCT}
Let $A$ be a separable \cstar-algebra and $G$ a metrizable compact group. Let $\alpha: G\to\Aut(A)$ be a strongly continuous action. Assume that $\alpha$ has the continuous Rokhlin property. Then $A$ $E$-dominates $A\rtimes_\alpha G$. In particular, the UCT passes from $A$ to its crossed product with respect to $\alpha$.
}
\begin{proof}
Since $\alpha$ has the continuous Rokhlin property, we can find an equivariant and unital $*$-homomorphism
\[
\phi: (\CC(G),\sigma) \to \bigl( F_{\fc,\alpha}(A), \tilde{\alpha}_\fc \bigl).
\]
Keeping in mind the observations \ref{central tensor} and \ref{path action}, this induces an equivariant and unital $*$-homomorphism
\[
\psi: (\CC(G)\otimes A,\sigma\otimes\alpha) \to \bigl( A_{\fc,\alpha}, \alpha_\fc \bigl)\quad\text{via}\quad f\otimes a \mapsto \phi(f)\cdot a.
\]
Consider the equivariant embedding $j: (A,\alpha)\to (\CC(G)\otimes A, \sigma\otimes\alpha)$ given by $a\mapsto\eins\otimes a$. Then obviously the composition $\psi\circ j$ coincides with the canonical embedding $\mu: A \into A_{\fc,\alpha}$. Since both of these $*$-homomorphisms are equivariant, this induces a commutative diagram
\[
\xymatrix{
A\rtimes_\alpha G \ar[rr]^{\hat{\mu}} \ar[dr]_{\hat{j}} && A_{\fc,\alpha}\rtimes_{\alpha_\fc} G \subset (A\rtimes_\alpha G)_\fc \\
& (\CC(G)\otimes A)\rtimes_{\sigma\otimes\alpha} G \ar[ur]_{\hat{\psi}} &
}
\]
By Green's imprimitivity theorem \cite[2.8]{Green80}, we have that $(\CC(G)\otimes A)\rtimes_{\sigma\otimes\alpha} G$ is isomorphic to $A\otimes\CK(L^2(G))$.

The above commutative diagram yields two $E$-theory elements $x=E(\hat{j})\in E(A\rtimes_\alpha G, A)$ and $y=E(\hat{\psi})\in E(A,A\rtimes_\alpha G)$ with $x\otimes y = E(\hat{\mu}) = 1$ in $E(A\rtimes_\alpha G, A\rtimes_\alpha G)$. In particular, we have verified \ref{E dominance} for $A$ and $A\rtimes_\alpha G$.
\end{proof}

{
\theorem \label{cRok fixed point UCT}
Let $A$ be a separable \cstar-algebra and $G$ a metrizable compact group. Let $\alpha: G\to\Aut(A)$ be a strongly continuous action. Assume that $\alpha$ has the continuous Rokhlin property. Then the canonical embedding $\iota: A^\alpha\to A$ has a right-inverse in $E$-theory.
In particular, $\iota$ induces $E$-dominance of $A$ over $A^\alpha$. Moreover, the UCT passes from $A$ to its fixed point algebra with respect to $\alpha$.
}
\begin{proof}
As we have observed in the proof of \ref{cRok crossed product UCT}, there exists an equivariant and unital $*$-homomorphism
\[
\psi: (\CC(G)\otimes A, \sigma\otimes\alpha) \to \bigl( A_{\fc,\alpha}, \alpha_\fc \bigl).
\]
Consider the group action $\sigma\otimes\alpha$ on $\CC(G)\otimes A \cong \CC(G,A)$, which is given (under this identification) by
\[
(\sigma\otimes\alpha)_g(f)(h) = \alpha_g(f(g^{-1}h))\quad\text{for all}~g,h\in G~\text{and}~f\in\CC(G,A).
\]
Now it is straightforward to check that a function $f\in\CC(G,A)$ is fixed under $\sigma\otimes\alpha$ if and only if it is of the form $f(h) = \alpha_h(a)$ for some $a\in A$. In particular, we have an isomorphism
\[
\kappa: A\to (\CC(G)\otimes A)^{\sigma\otimes\alpha}\quad\text{via}\quad \kappa(a)(h)=\alpha_h(a).
\]
Recall the equivariant embedding $j: A\to \CC(G)\otimes A$ given by $a\mapsto \eins\otimes a$.
From the definition of $\kappa$, it is obvious that one has $\kappa\circ\iota = j\circ\iota$. 

Observe that by compactness of $G$, the fixed point algebra of the path algebra $(A_{\fc,\alpha})^{\alpha_\fc}$ coincides with the path algebra of the fixed point algebra $(A^\alpha)_\fc$. Let $\mu: A^\alpha \into (A^\alpha)_\fc$ denote the canonical inclusion.

Combining all these observations, we obtain a commutative diagram
\[
\xymatrix{
A^\alpha \ar[rr]^\mu \ar[rd]_{\iota} && (A^\alpha)_\fc \\
& A \ar[ru]_{\psi\circ\kappa} &
}
\]
This yields two $E$-theory elements $x=E(\iota)\in E(A^\alpha, A)$ and $y=E(\psi\circ\kappa)\in E(A,A^\alpha)$ with $x\otimes y = E(\mu) = 1$ in $E(A^\alpha, A^\alpha)$. In particular, we have verified that $\iota$  induces $E$-dominance of $A$ over $A^\alpha$.
\end{proof}

%%%%%%%%%%%%%%%%%%%%%%%%%%%%%%%

\section{Circle actions and the continuous Rokhlin property}

In this section, we examine the special case of circle actions more closely. This is demonstrated in the next theorem:

{
\theorem \label{cRok circle}
Let $A$ be a separable, unital \cstar-algebra and $\alpha: \IT\to\Aut(A)$ a strongly continuous action. Assume that $\alpha$ has the continuous Rokhlin property. Then $A$ is $E$-equivalent to $A^\alpha\oplus SA^\alpha$.
}
\begin{proof}
Denote by $\iota: A^\alpha\to A$ the canonical embedding.
It follows from \cite[3.11]{Gardella14_1} that there is an automorphism $\theta\in\Aut(A^\alpha)$ such that there is an equivariant isomorphism
\[
\phi: (A^\alpha\rtimes_\theta\IZ, \hat{\theta}) \to (A,\alpha)\quad\text{with}\quad \phi|_{A^\alpha}=\iota.
\]
Identifying $A$ with $A^\alpha\rtimes_\theta\IZ$ by the above isomorphism, we apply the Pimsner-Voiculescu exact sequence in $KK$-theory \cite[19.6]{BlaKK}. For every separable \cstar-algebra $B$, there is a six-term exact sequence of the form
\[
\xymatrix{
KK(B,A^\alpha) \ar[rr]^{1-KK(\theta)} && KK(B,A^\alpha) \ar[rr]^{KK(\iota)} && KK(B,A) \ar[d] \\
KK(B,SA) \ar[u] && KK(B,SA^\alpha) \ar[ll]_{KK(S\iota)} && KK(B,SA^\alpha) \ar[ll]_{1-KK(S\theta)}
}
\]
Moreover, one can see from the construction of this sequence in \cite[Section 19]{BlaKK} that the vertical map on the right is given by right-multiplication with an element $z'\in KK(A,SA^\alpha)$, which is independant of $B$. The same is true for the left vertical map. Now $E$-theory naturally factors through $KK$-theory, in a way that is compatible with addition and forming Kasparov products, see \cite[25.5.8]{BlaKK}. Let $z$ denote the $E$-theory element induced by $z'$. Then we also have a Pimsner-Voiculescu exact sequence in $E$-theory:
\[
\xymatrix{
E(B,A^\alpha) \ar[rr]^{1-E(\theta)} && E(B,A^\alpha) \ar[rr]^{E(\iota)} && E(B,A) \ar[d]^{z} \\
E(B,SA) \ar[u] && E(B,SA^\alpha) \ar[ll]_{E(S\iota)} && E(B,SA^\alpha) \ar[ll]_{1-E(S\theta)}
}
\]

By \cite[3.11]{Gardella14_2}, it follows that $\theta$ is asymptotically representable. This particularly implies that $E(\theta)=1$ in $E(A^\alpha,A^\alpha)$.
Thus we can extract a short exact sequence
\[
\xymatrix{
0 \ar[r] & E(B,A^\alpha) \ar[r]^{E(\iota)} & E(B,A) \ar[r]^{z} & E(B,SA^\alpha) \ar[r] & 0.
}
\]
By \ref{cRok fixed point UCT}, there is some $y\in E(A,A^\alpha)$ with $E(\iota)\otimes y=1$ in $E(A^\alpha,A^\alpha)$. It follows that the map
\[
E(B, A)\to E(B, A^\alpha\oplus SA^\alpha),\quad x\mapsto x\otimes(y\oplus z)
\]
is an isomorphism of abelian groups for all $B$. Using surjectivity of this map in the case of $B=A^\alpha\oplus SA^\alpha$, we find $g\in E(A^\alpha\oplus SA^\alpha,A)$ with $g\otimes(y\oplus z)=1$ in $E(A^\alpha\oplus SA^\alpha, A^\alpha\oplus SA^\alpha)$. On the other hand, we have
\[
\begin{array}{ccl}
((y\oplus z)\otimes g)\otimes(y\oplus z) &=& (y\oplus z)\otimes (g\otimes(y\oplus z)) \\\\
&=& (y\oplus z)\otimes 1 = y\oplus z \\\\
&=& 1\otimes (y\oplus z).
\end{array}
\]
Hence, injectivity of the homomorphism 
\[
E(A,A)\to E(A,A^\alpha\oplus SA^\alpha),\quad x\mapsto x\otimes(y\oplus z)
\]
implies that $(y\otimes z)\otimes g=1$ in $E(A,A)$. This shows that $y\oplus z\in E(A,A^\alpha\oplus SA^\alpha)$ is an $E$-equivalence.
\end{proof}

\rem Applying the ordinary (i.e. $K$-theoretic) Pimsner-Voiculescu exact sequence as in the proof of \ref{cRok circle}, one recovers the isomorphisms $K_i(A)\cong K_i(A^\alpha)\oplus K_{i+1}(A^\alpha)$ for $i=0,1$. Moreover, in the case $i=0$, one can see that this isomorphism carries $[\eins_A]_0\in K_0(A)$ to the pair $([\eins_{A^\alpha}]_0, 0)\in K_0(A^\alpha)\oplus K_1(A^\alpha)$. 

\rem We further remark that unitality of $A$ is probably not necessary to obtain the result in \ref{cRok circle}. The only obstacle is to define a suitable notion of asymptotic representability for automorphisms on non-unital \cstar-algebras, and extending the results \cite[3.11]{Gardella14_1} and \cite[3.11]{Gardella14_2} to the non-unital case. Considering that this has been done in the finite group case (see \cite{Nawata13}), this should also be possible for circle actions in a similar fashion.

%\newpage
%%%%%%%%%%%%%%%%%%%%%%%%%%%%%%%%%%%%%%%%%%%%%%%%%%%%%%%%%%%%%%%%%%

\bibliographystyle{gabor}
\bibliography{master}

\end{document}